\documentclass[12 pt]{iopart}
\usepackage{iopams}
\usepackage{setstack}
\usepackage{latexsym}
\begin{document}
\eqnobysec
\title{An~Averaging~Theorem~for~Perturbed~KdV~Equation}
\author{HUANG Guan}
\address{C.M.L.S, Ecole Polytechnique, Palaiseau, France}
\ead{guan@math.polytechnique.fr}

\begin{abstract}
     We consider a perturbed KdV equation:
     \[\dot{u}+u_{xxx}-6uu_x=\epsilon f(x,u(\cdot)),\quad x\in \mathbb{T},\quad\int_{\mathbb{T}} u dx=0.\]
     For any periodic function $u(x)$, let $I(u)=(I_1(u),I_2(u),\cdots)\in\mathbb{R}_+^{\infty}$ be the vector, formed by the KdV integrals of motion, calculated for the potential $u(x)$. Assuming that the perturbation $\epsilon f(x,u(\cdot))$ is a smoothing mapping (e.g. it is a smooth function $\epsilon f(x)$, independent from $u$), and that solutions of the perturbed equation satisfy some mild a-priori assumptions, we prove that for solutions  $u(t,x)$  with typical initial data and  for $0\leqslant t\lesssim \epsilon^{-1}$,  the vector $I(u(t))$ may  be well approximated by a solution of the averaged equation.
\end{abstract}
\maketitle
\bibliographystyle{unsrt}
\setcounter{section}{-1}
\section{Introduction}
We consider a perturbed Korteweg-de Vries (KdV) equation  with zero mean-value  periodic boundary condition:
\begin{equation}
\label{pkdv}
\dot{u}+u_{xxx}-6uu_x=\epsilon f(x,u(\cdot)),\quad x\in \mathbb{T}=\mathbb{R}/\mathbb{Z},\quad\int_{\mathbb{T}}u(x,t)dx=0.
\end{equation}
 Here $\epsilon f(x,u(\cdot))$ is a nonlinear perturbation, specified below. For any $p\in \mathbb{R}$ we denote by $H^p$ the Sobolev space of order $p$, formed by  real-valued periodic functions with zero mean-value, provided with the homogeneous norm $||\cdot||_p$. Particularly,  if $p\in\mathbb{N}$ we  have
\[H^p=\Big\{u\in L^2(\mathbb{T}):||u||_p<\infty,\int_{\mathbb{T}}udx=0\Big\},\qquad
||u||_p^2=\int_{\mathbb{T}}\Big|\frac{\partial^p u}{\partial x^p}\Big|^2dx.\]
For any $p$, the operator $\frac{\partial}{\partial x}$ defines a linear isomorphism: $\frac{\partial}{\partial x}:$ $H^p\rightarrow H^{p-1}$. Denoting by $(\frac{\partial}{\partial x})^{-1}$  its inverse, we provide the spaces $H^p$, $p\geqslant 0$, with  a symplectic structure by means of the 2-form $\Omega$:
   \begin{equation}
   \Omega(u_1,u_2)=-\Big\langle (\frac{\partial}{\partial x})^{-1} u_1,u_2\Big\rangle,\label{symform1}
   \end{equation}
   where $\langle \cdot,\cdot\rangle$ is the scalar product in $L^2(\mathbb{T})$. Then in any space $H^p$, $p\geqslant1$,  the KdV equation (\ref{pkdv})$_{\epsilon=0}$   may be written as a Hamiltonian system with the Hamiltonian $\mathcal{H}$, given by
  $\mathcal{H}(u)=\int_{\mathbb{T}}\Big(\frac{1}{2}u_x^2+u^3\Big)dx$.
   That is, KdV   may be written as 
  \[\dot{u}=\frac{\partial}{\partial x}\nabla\mathcal{ H}(u).\]
  It is well-known that    KdV  is integrable. It means that the function space $H^p$ admits analytic symplectic coordinates $ v=(\mathbf{v}_1,\mathbf{v}_2,\cdots)=\Psi(u(\cdot))$, where $\mathbf{v}_j=(v_j,v_{-j})\in \mathbb{R}^2$, such that the quantities $I_j=\frac{1}{2}|\mathbf{v}_j|^2$, $j\geqslant1$, are actions (integrals of motion), while $\varphi_j=\mathrm{Arg} \; \mathbf{v}_j$, $j\geqslant1$, are angles. In the $(I,\varphi)$-variables,  KdV  takes the integrable form
\begin{equation}\dot{I}=0,\quad \dot{\varphi}=W(I),\label{eq0.2}
\end{equation}
where $W(I)\in \mathbb{R}^{\infty}$ is the frequency vector (see \cite{3,5}). The integrating  transformation  $\Psi$, called the  nonlinear Fourier transform,   for any $p\geqslant 0$  defines an analytic isomorphism $\Psi : H^p\rightarrow h^p$, where \[h^p=\Big\{ v=(\mathbf{v}_1,\mathbf{v}_2,\cdots): |v|_p^2=\sum_{j=1}^{+\infty}(2\pi j)^{2p+1}|\mathbf{v}_j|^2<\infty,\;\mathbf{v}_j\in \mathbb{R}^2,\; j\in\mathbb{N}\Big\}.\]

It is well established that  for a perturbed integrable finite-dimensional  system,
\[
\dot{I}=\epsilon f(I,\varphi), \quad \dot{\varphi}=W(I) + \epsilon g(I,\varphi),\quad \epsilon<<1,
\]
where $I\in \mathbb{R}^n$, $\varphi\in \mathbb{T}^n$, on  time intervals of order $\epsilon^{-1}$  the actions $I(t)$ may be well approximated by solutions of the averaged equation:
\[
\dot{J}=\epsilon \langle f\rangle(J),\quad \langle f\rangle(J)=\int_{\mathbb{T}^n} f(J,\varphi)d\varphi, 
\]
provided that the initial data $(I(0),\varphi(0))$ are typical (see \cite{4, 9, 10, 11}). This assertion is known as the \emph{ averaging principle}. But  
in the  infinite dimensional case, there is no similar  general result.    In \cite{1, 2},  S.~Kuksin and A.~Piatniski proved that the   averaging principle holds for the randomly perturbed KdV equation of the form:
\begin{equation}\dot{u}-\epsilon u_{xx} +u_{xxx}-6uu_x=\sqrt{\epsilon}\eta(t,x),\;\; x\in \mathbb{S}^1,\;\; \int u dx=\int \eta dx=0,\label{rmpkdv}
\end{equation}
where the  force $\eta$ is a white noise  in $t$,  is smooth in $x$ and is non-degenerate.  Our goal in this work is to justify the averaging principle for the KdV equation with deterministic perturbations, using the Anosov scheme (see \cite{4}), exploited earlier in the finite dimensional situation. The main technical difficulty to achieve this goal comes from the fact that to perform the scheme one has to use a measure in the function space which is quasi-invariant under the flow of the perturbed equation (it is needed to guarantee that a small 'bad' set which we have to prohibit for a solution of the perturbed equation at a time $t>0$ corresponds to a small set of initial data). For a reason, explained in Section~3, to construct such a quasi-invariant measure we have to assume that the perturbation $\epsilon f$ is smoothing. More  precisely, we assume that:
\smallskip

\noindent{\bf Assumption A.}  (i) For any $p\geqslant 0$,  the mapping defined by the perturbation in (\ref{pkdv}):
\begin{equation}
\mathcal{P}: H^{p}\rightarrow H^{p+\zeta_0},\quad u\mapsto f(x,u(\cdot)),
\end{equation}
is analytic. Here  $\zeta_0>1$ is a constant. 

(ii)  For any $p\geqslant 3$ and $T>0$, the perturbed KdV equation (\ref{pkdv})  with initial data
\[u(0)=u_0\in H^p,\]
 has a unique solution $u(t,x)\in H^p$ in the time interval $[-T\epsilon^{-1},T\epsilon^{-1}]$, and \[||u(t)||_p\leqslant C(p,||u_0||_p,T),\quad |t|\leqslant T\epsilon^{-1}.\]
\smallskip

We are mainly concerned with the behavior of the actions $I(u(t))\in \mathbb{R}_+^{\infty}$  for $|t|\lesssim\epsilon^{-1}$. For this end, it is convenient to pass to the  slow time $\tau=\epsilon t$ and    write the perturbed KdV equation (\ref{pkdv}) in the action-angle coordinates $(I,\varphi)$:
\begin{equation}
\frac{dI}{d\tau}= F(I,\varphi),\quad \frac{d\varphi}{d\tau}=\epsilon^{-1}W(I)+ G(I,\varphi).
\end{equation}
Here $I\in \mathbb{R}^{\infty}$, $\varphi\in \mathbb{T}^{\infty}$ and $\mathbb{T}^{\infty}:=\{\theta=(\theta_i)_{i\geqslant 1}, \theta_i\in \mathbb{T}\}$ is the infinite-dimensional torus, endowed with  the Tikhonov topology. The two functions  $F(I,\varphi)$ and $G(I,\varphi)$ are the perturbation term $\epsilon f$,   written in action-angle variables, see below (\ref{i1.3}) and (\ref{angle1}). The corresponding averaged equation is 
\begin{equation}
\frac{dJ}{d\tau}=\langle F\rangle(J),\quad \langle F\rangle(J)=\int_{\mathbb{T}^{\infty}}F(J,\varphi)d\varphi,\label{apkdv1}
\end{equation}
where $d\varphi$ is the Haar measure on $\mathbb{T}^{\infty}$.  It turns out that the (\ref{apkdv1}) is a Lipschitz equation, see below (\ref{l1}). 
We  denote by $h_{I+}^p$ the image of the  space $h^p$  under the action-mapping 
\[\pi_I: v\mapsto I,\quad I_j(v)=\frac{1}{2}|\mathbf{v}_j|^2,\quad j\geqslant 1.\]
Clearly, $I=\pi_I(v)\in h^p_{I+}\subset h^p_I$, where $h^p_I$ is the weighted $l^1$-space
\[h^p_I=\Big\{I\in \mathbb{R}^{\infty}:|I|_{h^p_I}=|I|_p=2\sum_{j=1}^{\infty} (2\pi j)^{2p+1}|I_j|< \infty\Big\}.\]
 and $h^p_{I+}$ is its positive octant, $h^p_{I+}=\{I\in h^p_I: I_j\geqslant 0, \;\forall j\}$. This is a closed subset of~$h^p_I$.
 
 For any $\theta=(\theta_i)_{i\geqslant1}\in \mathbb{T}^{\infty}$, let us denote by $\Phi_{\theta}$ the linear operator on the space of sequences $(\mathbf{v}_1, \mathbf{v}_2,\cdots)\in h^p$ which rotates each component $\mathbf{v}_j\in\mathbb{R}^2$ by the angle $\theta_j$.
 \smallskip
 
 \noindent{\bf Definition 0.1} A {\it Gaussian measure}  $\mu$ on the Hilbert space $h^p$ is said to be $\zeta_0$~-{\it admissible} (where $\zeta_0>1$ is the same as in assumption A),  if the following conditions are fulfilled:
 
 (i) It is non-degenerate and has zero mean value.
 
 (ii) It has a diagonal correlation operator $(\mathbf{v}_1,\mathbf{v}_2,\cdots)\mapsto (\sigma_1\mathbf{v}_1,\sigma_2\mathbf{v}_2,\cdots)$,  where every $\sigma_j>0$, $\sum_{j\geqslant1}\sigma_j<\infty$  and $j^{-\zeta_0}/\sigma_j=O(1)$. In particular,  $\mu$   is invariant under the rotations $\Phi_{\theta}$.
 \smallskip 
 
 Such measures can be written as:
 \begin{equation}\prod_{j=1}^{+\infty}\frac{(2\pi j)^{1+2p}}{2\pi\sigma_j}\exp\{-\frac{(2\pi j)^{1+2p}|\mathbf{v}_j|^2}{2\sigma_j}\}d\mathbf{v}_j,\label{mu0.10}
 \end{equation} where $d\mathbf{v}_j$, $j\geqslant 1$, is the Lebesgue measure on $\mathbb{R}^2$ (see \cite{13, 16}).  Clearly,  they are invariant under the KdV flow (\ref{eq0.2}).
 
 The main result of this work is the following theorem:
\smallskip

\noindent{\bf Theorem 0.2.} Fix any $p\geqslant 3$ and  $\bar{T}>0$. Let the curve $u^{\epsilon}(t)\in H^p$, $|t|\leqslant\epsilon^{-1}\bar{T}$ be a solution of equation (\ref{pkdv}) and $v^{\epsilon}(\tau)=\Psi(u^{\epsilon}(\epsilon^{-1}\tau))$, $\tau=\epsilon t$, $|\tau|\leqslant\bar{T}$. If assumption A is fulfilled and $\mu$ is a $\zeta_0$-admissible Gaussian measure on $h^p$, then

(i) For any $\rho>0$, there exists a Borel subset $\Gamma_{\rho}^{\epsilon}$ of $h^p$  and $\epsilon_{\rho}>0$ such that  $\lim_{\epsilon\rightarrow 0}\mu(h^p\setminus\Gamma^{\epsilon}_{\rho})=0$,  
and for $\epsilon\leqslant \epsilon_{\rho}$ we have
\begin{equation}
|I(v^{\epsilon}(\tau))-J(\tau)|_p\leqslant\rho,\quad
\mbox{for}\quad |\tau|\leqslant \bar{T},\quad v^{\epsilon}(0)\in\Gamma^{\epsilon}_{\rho},\label{main1}
\end{equation}
where $J(\tau)$, $|\tau|\leqslant \bar{T}$,  is a solution of the averaged equation (\ref{apkdv1}) with the inital data $J(0)=\pi_I(v^{\epsilon}(0))$.

(ii) There is a full measure subset $\Gamma_{\varphi}$ of $h^p$ with the following property:

\noindent If   $v^{\epsilon}(0)\in \Gamma_{\varphi}$, then for any $0\leqslant \bar{T}_1<\bar{T}_2\leqslant \bar{T}$   the image $\mu_{\bar{T}_1,\bar{T}_2}^{\epsilon}$ of the probability measure $(\bar{T}_2-\bar{T}_1)^{-1}d\tau$ on $[\bar{T}_1,\bar{T}_2]$ under the mapping   $ \tau\mapsto~\varphi(v^{\epsilon}(\tau))\in~\mathbb{T}^{\infty}$  converges weakly, as $\epsilon\rightarrow 0$, to the Haar measure $d\varphi$ on $\mathbb{T}^{\infty}$. 
 \medskip
 
  The assertion (ii) of the theorem means that  for any bounded continuous  function $g(\varphi)$ on $\mathbb{T}^{\infty}$,
\[\frac{1}{\bar{T}_2-\bar{T}_1}\int_{\bar{T}_1}^{\bar{T}_2}g(\varphi(v^{\epsilon}(\tau)))d\tau\rightarrow \int_{\mathbb{T}^{\infty}}g(\varphi)d\varphi,\quad \epsilon\rightarrow 0.\]

In particular, we have
\smallskip

{\noindent\bf Proposition 0.3.} The assumption A holds if in (\ref{pkdv}) $f=f(x) $ is a smooth function, independent from $u$.  
\smallskip

It is unknown for us that if the result of Theorem 0.2 remains true for equation (\ref{pkdv}) with non-smoothing perturbations, e.g. if the right hand side of equation (\ref{pkdv}) is $\epsilon u_{xx}$ or $-\epsilon u$. So we do not know whether a suitable analogy of the result in \cite{1,2} holds true if in equation (\ref{rmpkdv}) the noise $\eta$ vanishes.

The paper has the following structure: Section 1 is about the  transformation which integrates the KdV and its Birkhoff normal form. In Section 2 we discuss the averaged equation. We prove that the $\zeta_0$-admissible Gaussian measures are quasi-invariant under the flow of equation (\ref{pkdv}) in Section 3. Finally in Section 4  and Section 5 we establish  the main theorem and Proposition ~0.3.

\noindent{\bf Agreements.} Analyticity of maps $B_1\rightarrow B_2$ between Banach spaces $B_1$ and $B_2$, which are the real parts of complex spaces $B_1^c$ and $B_2^c$, is understood in the sense of Fr\'echet. All analytic maps that we consider possess the following additional property: for any $R$, a map  extends to a bounded analytical mapping in a complex ($\delta_R>0$)-neighborhood of the ball $\{|u|_{B_1}<R\}$ in $B_1^c$.

\noindent{\bf Notation.} We use   capital letters  $C$ or $C(a_1,a_2,\dots)$ to denote positive constants  that depend on the parameters~$a_1$,~$a_2$, $\dots$ but not on the unknown function $u$. We denote $\mathbb{Z}_{\geqslant0}=\{n\in\mathbb{Z},n\geqslant 0\}.$ For an infinite-dimensional vector $w=(w_1,w_2,\dots)$ and any $n\in\mathbb{N}$ we denote $w^n=(w_1,\dots,w_n,0,0,\dots)$. We often identify $w^n$ with a corresponding $n$-vector.

 \section{Preliminaries on the KdV equation}
     In this section we discuss integrability of the KdV equation (0.1)$_{\epsilon=0}$.

     \subsection{ Nonlinear Fourier transform for KdV}
We provide the $L^2$-space $H^0$ with the Hilbert basis $\{ e_s, s\in \mathbb{Z}\setminus\{0\}\}$,
     \[
     e_s =\cases{\sqrt2\cos(2\pi s x) \quad s>0,\\
     \sqrt2\sin(2\pi sx)\quad s<0. }
     \]

\noindent {\bf Theorem 1.1.}  There exists an analytic diffeomorphism $\Psi:H^0 \mapsto h^0$ and an analytic functional $K$ on $h^0$ of the form $K(v)=\tilde{K}(I(v))$, where the function $\tilde{K}(I)$ is analytic in a suitable neighborhood of the octant $h_{I+}^0$ in  $h_I^0$, with the following properties:
 
 (i) The mapping $\Psi $ defines an analytic diffeomorphism $\Psi : H^p \mapsto h^p$, for any $p\in \mathbb{Z}_{\geqslant0}$. This is a symplectomorphism of the spaces $(H^p,\Omega)$ (see (\ref{symform1})  and $(h^p,\omega_2)$, where  $ \omega _2=\sum dv_k\wedge dv_{-k}$.
 
 (ii) The differential  $d\Psi(0)$ takes the form $\sum u_se_s\mapsto v,v_s=|2\pi s|^{-1/2}u_s$.
 
(iii) A curve $u\in C^1(0,T;H^0)$ is a solution of the KdV equation (\ref{pkdv})$_{\epsilon=0}$ if and only if $v(t)=\Psi(u(t))$ satisfies the equation
  
\begin{equation} \dot{\mathbf{v}}_j=\left(\begin{array}{cc} 0 & -1\\1& 0 \end{array}\right)\frac{\partial \tilde{K}}{\partial I_j}(I)\mathbf{v}_j  , \quad\mathbf{v}_j=(v_j,v_{-j})\in\mathbb{R}^2,\;j\in \mathbb{N}.
 \end{equation}
 \smallskip
 
Since the maps $\Psi$ and $\Psi^{-1}$ are analytic, then for $m=0,1,2\dots $, we have
\[||d^j\Psi(u)||_m\leqslant P_m(||u||_m),\quad
 ||d^j\Psi^{-1}(v)||_m\leqslant Q_m(|v|_m ), \quad j=0,1,2,
\]
where $P_m$ and $Q_m$ are continuous functions (cf.  the agreements).

 We denote  \[W(I)=(W_1, W_2,\dots),\; \; W_k(I)=\frac{\partial \tilde{K}}{\partial I_k}(I),\;\;   k=1,2,\dots .\] 

\noindent {\bf Lemma 1.2.} For any $n\in \mathbb{	N}$, if $ I_{n+1}=I_{n+2}=\cdots=0$, then
\[\mbox{det}\Big((\frac{\partial W_i}{\partial I_j})_{1\leqslant i,j\leqslant n}\Big)\neq 0.\]
Let $l^{\infty}_{-1}$ be the Banach space of all real sequences $l=(l_1,l_2,\dots)$ with the norm
 \[|l|_{-1}=\sup_{n\geqslant 1}n^{-1}|l_n|<\infty.\]
 Denote $\kappa=(\kappa_n)_{n\geqslant 1}$, where
 $\kappa_n=(2\pi n)^3$.
 \smallskip
 
 \noindent{\bf Lemma 1.3.} The normalized frequency map
\[\tilde{W}: I\mapsto \tilde {W}(I)=W(I)-\kappa \]
 is real analytic as a map from $h^1_{I+}$ to $l^{\infty}_{-1}$.  
 \medskip
 
 The coordinates $ v=\Psi(u) $ are called the \emph{Birkhoff coordinates}, and the form (1.1) of KdV is its \emph{Birkhoff normal form}.
 See \cite{3} for Theorem 1.1 and Lemma 1.3. A detailed proof of 
 Lemma 1.2 can be found in \cite{5}.

 \subsection{ Equation (\ref{pkdv}) in the Birkhoff coordinates.}
 For $k=1,2\dots $  we denote:\[\Psi_k : H^m \to \mathbb{R}^2,\quad \Psi_k(u)=\mathbf{v}_k,\]
 where  $\Psi(u)=v=(\mathbf{v}_1,\mathbf{v}_2,\dots)$.
  Let $u(t)$ be a solution of equation (\ref{pkdv}). We get
 \begin{equation}
 \dot{\mathbf{v}}_k=d\Psi_k(u)(\epsilon f(x,u)+V(u)), \quad k\geqslant 1,\label{tvpkdv}
 \end{equation}
 where $V(u)=-u_{xxx}+6uu_x$. Since $I_k(v)=\frac{1}{2}|\Psi_k|^2$ is an integral of motion of KdV equation (\ref{pkdv})$_{\epsilon=0}$, we have 
 \begin{equation}\dot{I}_k=\epsilon (d\Psi_k(u)f(x,u),\mathbf{v}_k):=\epsilon F_k(v).\label{i1.3}
 \end{equation}
Here and below $(\cdot,\cdot)$ indicates the scalar product in $\mathbb{R}^2$.

For $k\geqslant 1$ define $\varphi_k=\arctan(\frac{v_{-k}}{v_k})$  if $\mathbf{v}_k\neq 0$, and $\varphi_k=0$ if $\mathbf{v}_k=0$. Using equation (1.1), we get
\begin{equation}\dot{\varphi}_k=W_k(I)+\epsilon|\mathbf{v}_k|^{-2}(d\Psi_k(u)f(x,u),\mathbf{v}_k^{\bot}),\quad \mbox{if} \quad\mathbf{v}_k\neq 0,\label{angle1}
\end{equation}
  where $\mathbf{v}_k^{\bot}=(-v_{-k},v_k) $. Denoting for brevity, the vector field in equation (1.4) by $W_k(I)+\epsilon G_k(v)$, we rewrite the equation for the pair $(I_k, \varphi_k) (k\geqslant1)$ as
 \begin{eqnarray}
 \eqalign{\dot{I}_k(t)=\epsilon F_k(v) =\epsilon F_k(I,\varphi), \\ 
 \dot{\varphi}_k(t)=W_k(I)+\epsilon G_k(v). }
 \label{iapkdv}
 \end{eqnarray}
 
  We set \[F(I,\varphi)=(F_1(I,\varphi),F_2(I,\varphi),\dots).\]

  In the following lemma $P_k$ and $P_k^j$ are some fixed continuous functions.
 \smallskip
 
 \noindent{\bf Lemma 1.4.} For $k,  j\in \mathbb{N}$,   we have for any $p\geqslant 0$
 \begin{enumerate}
 \item[(i)] The function $F_k(v)$ is analytic in each space $h^p$.
 \item[(ii)] For any $p\geqslant 0$, $\delta>0$, the function $G_k(v)\chi_{ \{I_k \geqslant \delta\}}$ is bounded by $\delta^{-1/2}P_k(|v|_p)$.
 \item[(iii)] For any $\delta>0$, the function $\frac{\partial F_k}{\partial I_j}(I,\varphi)\chi_{\{I_j\geqslant \delta\} }$ is bounded by $\delta^{-1/2}P_k^j(|v|_p)$.
 \item[(iv)] The function $\frac{\partial F_k}{\partial \varphi_j}(I,\varphi)$ is bounded by $P_k^j(|v|_p)$,  and  for any $n\in \mathbb{N}$ and  $(I_1,\dots,I_n)\in \mathbb{R}^n_+$, the function  $F_k(I_1,\varphi_1,\dots, I_n,\varphi_n,0,\dots)$ is analytic on $\mathbb{T}^n$.
 
  \end{enumerate}
  
  \noindent{\it Proof}:  Items (i) and (ii) follow directly from Theorem 1.1. Items  (iii) and (iv) follow  from item (i) and the chain-rule:
  \[\eqalign{\frac{\partial F_k}{\partial \varphi_j}=\sqrt{2I_j}\Big(\frac{\partial F_k}{\partial v_{-j}}\cos(\varphi_j)-\frac{\partial F_k}{\partial v_{j}}\sin(\varphi_j)\Big),\\
  \frac{\partial F_k}{\partial I_j}=(\sqrt{2I_j})^{-1}\Big(\frac{\partial F_k}{\partial v_j}\cos(\varphi_j)+\frac{\partial F_k}{\partial v_{-j}}\sin(\varphi_j)\Big).\quad \square}
  \]
  
  From this lemma we know that equation (\ref{iapkdv}) may have singularities at $\partial h_{I+}^p$.
  We  denote
   \[\eqalign{\Pi_I: h^p\to h^p_I,\quad \Pi_I(v)=I(v),\\
   \Pi_{I,\varphi}:h^p\to h_I^p\times \mathbb{T}^{\infty},\quad \Pi_{I,\varphi}(v)=(I(v),\varphi(v)).}
   \]
   Abusing notation, we will identify $v$ with $(I,\varphi)=\Pi_{I,\varphi}(v)$.
   \smallskip
   
   \noindent{\bf Definition 1.5.} For $p\geqslant 3$,  we say that a curve   $(I(t),\varphi(t))$, $|t|\leqslant T $, is a regular solution of equation (\ref{iapkdv}),   if there  exists a solution $u(t)\in H^p$ of equation (\ref{pkdv})  such that 
   $u(t)\in H^p$ and \[\Pi_{I,\varphi}(\Psi(u(t)))=(I(t),\varphi(t)),\quad |t|\leqslant T.\]
  \smallskip

   If $(I(t),\varphi(t))$ is a regular solution of (\ref{iapkdv}) and $|I(0)|_p\leqslant M_0$, then by assumption A we have\begin{equation}|I(t)|_p=|v(t)|_p^2\leqslant C(p,M_0,T),\quad |t|\leqslant T\epsilon^{-1}.\label{le 1.6}
   \end{equation}

\section{ Averaged equation }

     For a function $f$ on a Hilbert space $H$, we write $f\in Lip_{loc}(H)$ if
      \begin{equation} |f(u_1)-f(u_2)|\leqslant P(R)||u_1-u_2||, \quad  ||u_1||, ||u_2||\leqslant R, \label{2.1}
      \end{equation}
    for a suitable continuous function $P$ which depends on $f$.
      Clearly, the set of functions $Lip_{loc}(H)$ is an algebra. By  the Cauchy inequality, any analytic function on $H$ belongs to $Lip_{loc}(H)$ (see agreements).  In particularly, for any $k\geqslant 1$,
     \[\eqalign{W_k(I)\in Lip_{loc}(h_I^p),\;\;   p\geqslant 1,\quad \mbox{and} \quad
      F_k(v)\in Lip_{loc}(h^p), \;\; p\geqslant 0.}
      \]
      
     In the further analysis, we systematically use the fact that the functional $F_k(v)$ only weakly depends on the tail of the  vector $v$. Now we state the corresponding  results.
      Let $f\in Lip_{loc}(h^{p})$ and $v\in h^{p_1}$,   $ p_1>p$. Denoting by $\Pi^M$,  $M\geqslant 1$ the projection
     \[ \Pi^M: h^0 \rightarrow h^0, \quad (\mathbf{v}_1,\mathbf{v}_2,\dots) \mapsto (\mathbf{v}_1,\dots,\mathbf{v}_M, 0,\dots),\]
     we have $|v-\Pi^M v|_{p}\leqslant (2\pi M)^{-(p_1-p)}|v|_{p_1}$.
     Accordingly, 
     \begin{equation}
     |f(v)-f(\Pi^M v)|\leqslant P(|v|_{p_1})(2\pi M)^{-(p_1-p)}.\label{lip}
     \end{equation}
     
     The torus $\mathbb{T}^M$ acts on the space $\Pi_M h^0$ by linear transformations $\Phi_{\theta_M}$, $\theta_M \in \mathbb{T}^M$, where $\Phi_{\theta_M}: (I_M, \varphi_M)\mapsto (I_M, \varphi_M+\theta_M)$. Similarly, the torus $\mathbb{T}^{\infty}$ acts on $h^0$ by linear transformations $\Phi_{\theta}: (I, \varphi) \mapsto(I,\varphi +\theta)$ with $\theta\in \mathbb{T}^{\infty}$.
     
     For a function $f\in Lip_{loc}(h^{p})$ and a positive integer $N$  we define the average of $f$ in the first $N$ angles as the function
     \[\langle f\rangle_N(v)=\int_{\mathbb{T}^N} f\Big((\Phi_{\theta_N}\oplus \mbox{Id})(v)\Big)d\theta_N,\]
     and define the averaging in all angles as \[\langle f\rangle(v)=\int_{\mathbb{T}^{\infty}}f(\Phi_{\theta}(v))d\theta,\]
     where  $d\theta$ is the Haar measure on $\mathbb{T}^{\infty}$.  The estimate (\ref{lip}) readily implies that \[|\langle f\rangle_N(v)-\langle f\rangle(v)|\leqslant P(R) 	(2\pi N)^{-(p_1-p)},  \quad |v|_{p_1}\leqslant R.\]
     Let $v=(I, \varphi)$, then $\langle f\rangle_N $ is a function independent of $\varphi_1, \cdots, \varphi_N$, and $\langle f\rangle$ is independent of $\varphi$. Thus $\langle f\rangle$ can be written as $\langle f\rangle(I)$.
     \smallskip
     
     \noindent{\bf Lemma 2.1.} (See \cite{1}). Let $f\in Lip_{loc}(h^{p})$, then
     \begin{enumerate}
     \item[(i)] The functions $\langle f\rangle_N(v)$ and $\langle f\rangle(v)$ satisfy (\ref{2.1}) with the same function $P$ as $f$ and take the same value at the origin.
     \item[(ii)] These two functions  are smooth (analytic) if $f$ is. If $f$ is smooth, then $ \langle f\rangle(I) $ is a smooth function with respect to vector $(I_1,\dots, I_M)$, for any $	M$. If $f(v)$ is analytic in the space $h^{p}$, then $\langle f\rangle(I)$ is analytic in the space $h_I^{p}$.
     \end{enumerate}
     \medskip

   We recall that a vector $\omega\in\mathbb{R}^n$ is  \emph{non-resonant} if
   \[\omega\cdot k\neq 0, \quad\forall k\in \mathbb{Z}^n\setminus\{0\}.\]
   Denote by $C^{0+1}(\mathbb{T}^n)$ the set of all Lipschitz functions on $\mathbb{T}^n$.
   \smallskip
   
   \noindent{\bf Lemma 2.2.} Let $f\in C^{0+1}(\mathbb{T}^n)$ for some $n\in\mathbb{N}$. Then for any non-resonant vector $\omega\in\mathbb{R}^{n}$  we have
   \[\lim_{T\rightarrow\infty}\frac{1}{T}\int_0^T f(x_0+\omega t)dt=\langle f\rangle,\]
   uniformly in $x_0\in\mathbb{T}^n$. The rate of convergence depends on $n$, $\omega$ and $f$.   
   
   \noindent{\it Proof}:\quad Let us write $f(x)$ as the Fourier series $f(x)=\sum f_ke^{ik\cdot x}$. Since the Fourier series of a Lipschitz function converges uniformly (see \cite{12}), for any $\epsilon>0$ we may find $R=R_{\epsilon}$ such that $\Big|\sum_{|k|>R}f_ke^{ik\cdot x}\Big|\leqslant\frac{\epsilon}{2}$ for all $x$. Now it is enough to show that
   \begin{equation}
   \Big|\frac{1}{T}\int_0^T f_R(x_0+\omega t)dt-f_0\Big|\leqslant \frac{\epsilon}{2},\quad \forall T\geqslant T_{\epsilon},
   \end{equation}
   for a suitable $T_{\epsilon}$, where $f_R(x)=\sum_{|k|\leqslant R}f_ke^{ik\cdot x}$. Observing that 
   \[\Big|\frac{1}{T}\int_0^T e^{ik\cdot(x_0+\omega t)}dt\Big|\leqslant \frac{2}{T|k\cdot \omega|},\]
   for each nonzero $k$. Therefore the l.h.s of (2.3) is smaller than
   \[\frac{2}{T}\Big(\inf_{|k|\leqslant R}|k\cdot \omega|\Big)^{-1}\sum_{|k|\leqslant R}|f_k|.\]
   
   The assertion of the lemma follows.\quad $\square$

\section{Quasi-invariance of Gaussian measures}

Fix any integer $p\geqslant 3$, and let $\mu$ be a  $\zeta_0$-admissible Gaussian measure on the Hilbert space $h^p$.  In this section we will discuss how this measure evolves under the flow of the perturbed KdV equation (\ref{pkdv}). We follow a classical procedure based on  finite dimensional approximations (see e.g.  \cite{14,16}).

We suppose the assumption A holds. Let us write the equation (\ref{pkdv}) in the Birkhoff normal form, using  the slow time $\tau=\epsilon t$:
\begin{equation}
\frac{d}{d\tau}\mathbf{v}_j=\epsilon^{-1}\mathcal{J}W_j(I)\mathbf{v}_j 
      +\mathbf{X}_j(v),\quad j\in\mathbb{N}, \label{vpkdv}
     \end{equation}
where $\mathbf{ X}_j=(X_j,\;X_{-j})^{t}\in \mathbb{R}^2$ and $\mathcal{J}=\left(\begin{array}{cc}0 & -1 \\0 & 0\end{array}\right)$. 

For any $n\in\mathbb{N}$, we consider the $2n$-dimensional subspace $\pi_n(h^p)$ of $h^p$ with coordinates $v^n=(\mathbf{v}_1,\dots,\mathbf{v}_n,0,\dots)$.  On $\pi_n(h^p)$, we define the following finite-dimensional systems:  
\begin{equation}
\frac{d}{d\tau}\vec{\omega}_j=\epsilon^{-1}\mathcal{J}W_j(I(\omega^n))\vec{\omega}_j
      +\mathbf{X}_j(\omega^n),\quad 1\leqslant j\leqslant n, \label{nvpkdv}
 \end{equation}
 where $\vec{\omega}_j=(\omega_j,\omega_{-j})^t\in\mathbb{R}^2$ and $\omega^n=(\vec{\omega}_1,\dots,\vec{\omega}_n,0,\dots)\in \pi_n(h^p)$.
 
 We denote  $X^n(v^n)=(\mathbf{X}_1(v^n),\dots, \mathbf{X}_n(v^n),0,\dots)$ and $X(v)=(\mathbf{X}_1(v),\dots)$.
By assumption A and Theorem 1.1,  for any $p\geqslant 0$ the mapping 
 \begin{equation}X:h^{p}\to h^{p+\zeta_0},\quad v\mapsto X(v)\;\mbox{ is analytic.}\label{prop 3.1}
 \end{equation}
\smallskip

\noindent{\bf Theorem 3.2.} For any $T>0$, $\omega^n(\cdot)$  converges to $v(\cdot)$ as $n\rightarrow \infty$ in $C([-T,T]; h^p)$, where $v(\cdot)$ and $\omega^n(\cdot)$ are,  respectively,  solutions of (\ref{vpkdv}) and (\ref{nvpkdv}) with initial data $v(0)\in h^p$ and $\omega^n(0)=v^n(0)\in \pi_n(h^p)$.

\noindent{\it Proof}: \quad
Fix any $M_0>0$. From (\ref{le 1.6}) we know that there exists  a constant $M_1$ such that if $|v(0)|_p\leqslant M_0$, then\begin{equation}
|v(\tau)|_p\leqslant M_1,\quad \tau\in[0,T].
\end{equation}
The equation (\ref{nvpkdv}) yields that 
\begin{equation}
\frac{d}{d \tau} |\omega^n|_p^2=2\sum_{j=1}^{n}j^{1+2p}\vec{\omega}_j \cdot \mathbf{X}_j(\omega^n):=\chi^n(\omega^n).\label{chi}
\end{equation}
We  define
\[
\chi(v):=2\sum_{j=1}^{\infty}j^{1+2p}\mathbf{v}_j\cdot \mathbf{X}_j(v).
\]
By (\ref{prop 3.1}), we know that  there exists  a constant $C_1>0$ such that
\begin{equation}
|\chi^n(\omega^n)|\leqslant C_1,  \quad |\omega^n|_p\leqslant 2M_1,\quad\forall n\in\mathbb{N}.
\end{equation}
Denote $\bar{\tau}=M_1/C_1$, then if  $|\omega^n(0)|_p\leqslant M_0$, then 
\begin{equation}
|\omega^n(\tau)|_p\leqslant 2M_1, \quad \tau\in[-\bar{\tau},\bar{\tau}], \quad \forall n\in\mathbb{N}.
\end{equation}

\noindent{\bf Lemma 3.3.}\quad In the space $C([-\bar{\tau},\bar{\tau}],h^{p-1})$,  we have the convergence 
\[\omega^n(\cdot)\rightarrow v(\cdot)\quad  \mbox{as}\quad n\rightarrow \infty.\]

\noindent{\it Proof}:\quad Denote $\vec{\xi}_j=\mathbf{v}_j-\vec{\omega}_j$, $I_v=I(v)$ and $I_{\omega^n}=I(\omega^n)$. Since $\mathcal{J}\mathbf{v}_j=\mathbf{v}_j^{\bot}$,  using  equations~(\ref{vpkdv}) and (\ref{nvpkdv}), for $1\leqslant j\leqslant n$, we get 
\[
\eqalign{
\frac{d}{d\tau}|\vec{\xi}_j|^2
=2(\vec{\xi}_j)^t[\epsilon^{-1}\mathcal{J}(W_j(I_v)\mathbf{v}_j-W_j(I_{\omega^n})\vec{\omega}_j)+\mathbf{X}_j(v)-\mathbf{X}_j(\omega^n)]\\
=2\epsilon^{-1}[W_j(I_v)-W_j(I_{\omega^n})]\mathbf{v}_j\cdot(\vec{\omega}_j)^{\bot}+2(\vec{\xi}_j)^t\cdot(\mathbf{X}_j(v)-\mathbf{X}_j(\omega^n)).
}
\]
By Lemma 1.3 and Cauchy's inequality, we know that
\[
\Big|W_j(I(v))-W_j(I(\omega^n))\Big|\leqslant C_2 (M_1)j|v-\omega^n|_{p-1}.
\]
Using (\ref{prop 3.1}) we get that
\[
\frac{d}{d\tau}|v-\omega^n|^2_{p-1}\leqslant C_3(\epsilon, M_1)|v-\omega^n|^2_{p-1}+a_n(v),\quad \tau\in[-\bar{\tau},\bar{\tau}],
\]
where 
\[
a_n(v)=\sum_{j=n+1}^{\infty}j^{2p-1}\mathbf{v}_j\cdot \mathbf{X}_j(v).
\]
Obviously, $a_n(v)\rightarrow 0$ as $n\rightarrow \infty$ uniformly for $|v|_p\leqslant M_1$.

The lemma now follows directly from Gronwall's Lemma.\quad $\square$
\smallskip

\noindent{\bf Lemma 3.4.} If $\omega^n(0)\rightarrow v(0)$ strongly in $h^p$ and $\tau_n\to\tau$, $\tau_n\in[-\bar{\tau},\bar{\tau}]$, as $n\rightarrow\infty$, then \[\lim_{n\rightarrow \infty}|v(\tau)-\omega^n(\tau_n)|_p= 0.\] 

\noindent{\it Proof}:\quad From (\ref{chi}) we know that for any $\tau_n\in [-\bar{\tau},\bar{\tau}]$, 
\[
|\omega^n(\tau_n)|^2_p-|\omega^n(0)|^2_p=\int_0^{\tau_n}\chi^n(\omega^n(s))ds.
\]
Since $\omega^n(0)\rightarrow v(0)$ strongly in $h^p$,  then using (\ref{prop 3.1}) and Lemma~3.3  we get
\[
\eqalign{
|v(\tau)|^2_p\leqslant \liminf_{n\rightarrow \infty}|\omega^n(\tau_n)|^2_p\leqslant \limsup_{n\rightarrow\infty}|\omega^n(\tau_n)|^2_p\\
=\limsup_{n\rightarrow\infty}\Big(|\omega^n(0)|_p^2+\int_0^{\tau_n}\chi^n(\omega^n(s))ds\Big)=|v(0)|_p^2+\int_0^{\tau}\chi(v(s))ds\\
=|v(\tau)|_p^2.
}
\]
Therefore,
$\lim_{n\rightarrow\infty}|\omega^n(\tau_n)|_p=|v(\tau)|_p.$  
Since $\omega^n(\tau_n)\rightarrow v(\tau)$ in the space $h^{p-1}$ as $n\rightarrow\infty$, then the required convergence follows.~\quad~$\square$
\smallskip

\noindent{\bf Lemma 3.5.} In the space $C([-\bar{\tau},\bar{\tau}],h^p)$, $\omega^n(\cdot)\rightarrow v(\cdot)$ as $n\rightarrow\infty$.

\noindent{\it Proof}: \quad Suppose this statement is invalid. Then there exists $\delta>0$ and a sequence $\{\tau^n\}_{n\in\mathbb{N}}\subset[-\bar{\tau},\bar{\tau}]$ such that 
\[|\omega^n(\tau^n)-v(\tau^n)|_p\geqslant\delta.\]
Let $\{\tau^{n_k}\}_{k\in\mathbb{N}}$ be a subsequence of the sequence $\{\tau^n\}_{n\in\mathbb{N}}$ converging to some $\tau^0\in [-\bar{\tau},\bar{\tau}]$. But $v(\tau^{n_k})\rightarrow v(\tau^0)$ in $h^p$ as $k\rightarrow\infty$, and using Lemma ~3.4, we can get $\omega^{n_k}(\tau^{n_k})\rightarrow v(\tau^0)$ as $k\rightarrow \infty$ in $h^p$. So we get a contradiction, and Lemma 3.5 is proved.\quad $\square$

If $T\leqslant \bar{\tau}$, the theorem is proved, otherwise we iterate the above procedure.
This finishes the  proof of Theorem 3.2. \quad $\square$ 
\medskip

Let  $\mathcal{S}_v^{\tau}$ denote the flow determined by  equations (\ref{vpkdv}) in the space  $h^p$, and 
\[B^v_p(M):=\{v\in h^p: |v|_p\leqslant M\}.\]

\noindent {\bf Theorem 3.6.} For any $M_0>0$ and  $T>0$, there exists a constant $C>0$ which depends only on $M_0$ and $T$, such that if 
$A$ is a open subset of $B_p^v(M_0)$, then for $\tau\in [0,T]$, we have 
\[e^{-C\tau}\mu(A)\leqslant \mu(\mathcal{S}_v^{\tau}(A))\leqslant e^{C\tau}\mu(A).\] 

\noindent{\it Proof}:\quad From (\ref{le 1.6}) we know that there is constant $M_1$ which  only depends on $M_0$ and $T$, such that if $v(0)\in B_p^v(M_0)$, then 
\begin{equation}
v(\tau)\in B^v_p(M_1), \quad |\tau|\leqslant T.\label{m1}
\end{equation}
For any $n\in\mathbb{N}$, consider the  measure $\mu_n=\pi_n\circ\mu$ on the subspace $\pi_n(h^p)$.
Since $\mu$ is a $\zeta_0$-admissible Gaussian measure,  by (\ref{mu0.10}) $\mu_n$ has the following density with respect to  the Lebesgue measure:
\[b_n(v^n):=(2\pi)^{-n}\prod_{j=1}^n(2\pi j)^{1+2p}\sigma_j^{-1}\exp\{-\frac{1}{2}\sum_{j=1}^n\frac{j^{1+2p}|\mathbf{v}_j|^2}{\sigma_j}\}.\]
Let $\mathcal{S}_n^{\tau}$ be the flow determined by equations (\ref{nvpkdv}) on subspace $\pi_n(h^p)$.
For any open set  $A_n\subset \pi_n(B^v_p(M_0))$, due to Theorem  A  in the appendix, we have 
\[
\eqalign{
\frac{d}{d\tau}\mu_n(\mathcal{S}^{\tau}_n(A_n))\\
=\int_{\mathcal{S}^{\tau}_n(A_n)}\sum_{j=1}^{n}\Big(\frac{\partial (b_n(v^n)X_j(v^n))}{\partial v_j}+\frac{\partial(b_n(v^n)X_{-j}(v^n))}{\partial v_{-j}}\Big)dv^n\\
=\int_{\mathcal{S}^{\tau}_n(A_n)}\sum_{j=1}^{n}j^{2p+1}\Big(\frac{v_jX_j+v_{-j}X_{-j}}{\sigma_j}+\frac{\partial X_j}{\partial v_j}+\frac{\partial X_{-j}}{\partial v_{-j}}\Big)b_n(v^n)dv^n\\
:=\int_{\mathcal{S}^{\tau}_n(A_n)}c^n(v^n)b_n(v^n)dv^n
}
\]
Since $j^{-\zeta_0}/\sigma_j=O(1)$, using  (\ref{prop 3.1}) and the Cauchy's inequality,  there exists a constant $C$ which depends only on $M_1$, such that 
\begin{equation}
|c^n(v^n)|\leqslant C,\quad v^n\in \pi_n(B_p^v(M_1),\quad \forall n\in\mathbb{N}.\label{uniformn}
\end{equation}
 We have
\begin{equation}
e^{-C\tau}\mu_n(A_n)\leqslant \mu_n(\mathcal{S}_n^{\tau}(A_n))\leqslant e^{C\tau}\mu_n(A_n),\label{m2}
\end{equation}
as long as $\mathcal{S}_n^{\tau}(A_n)\subset \pi_n(B_p^v(M_1))$.

Since $\mu_n$ convergences weakly to $\mu$,  the theorem follows from (\ref{m1}), (\ref{m2}) and Theorem 3.2 (see \cite{14, 16}).\quad$\square$

  \section{ Proof of the main theorem }
       In this section we prove Theorem 0.2 by developing a suitable infinite-dimensional version of the Anosov scheme (see \cite{4,9,10,11}), and by studying  the behavior of the regular solutions of equation (\ref{iapkdv}) and the corresponding solutions of (\ref{pkdv}). We fix $p\geqslant 3$.  Assume  $u(0)=u_0\in H^p$.
       So \begin{equation}
       \Pi_{I,\varphi}(\Psi(u_0))=(I_0,\varphi_0)\in h_{I+}^p\times \mathbb{T}^{\infty},\quad p\geqslant 3.
       \end{equation}   
       \subsection{Proof of the assertion (i)}    
       We denote \[B_p^I(M)=\{I\in h^p_{I+}: |I|_p\leqslant M\}.\]
      Without loss of generality, we  assume that $\bar{T}=1$ and $t\geqslant 0$.
                 
      Fix any  $M_0>0$.
      Let \[(I_0,\varphi_0)\in B^I_p(M_0)\times \mathbb{T}^{\infty}:=\Gamma_0,\] 
      that is,
      \[v_0=\Psi(u_0)\in B_p^v(\sqrt{M_0}).\]
      Let $(I(t),\varphi(t))$ be a regular solution of the system (\ref{iapkdv}) with $(I(0),\varphi(0))=(I_0,\varphi_0)$. Then by (\ref{le 1.6}), there exists $M_1\geqslant M_0$  such that
       \begin{equation}I(t)\in B_p^I(M_1),\quad t\in [0,\epsilon^{-1}].
       \end{equation}
      
           By the definition of the perturbation we know that 
           \begin{equation}
           |F(I,\varphi)|_1\leqslant C_{M_1},\quad \forall (I,\varphi)\in B_p^I(M_1)\times \mathbb{T}^{\infty},\label{boundf}
           \end{equation}
           where the  constant  $C_{M_1}$ depends only on $M_1$.

       We denote $I^m=(I_1,\dots, I_m, 0,0,\dots)$,  $\varphi^m=(\varphi_1,\dots,\varphi_m,0,0,\dots)$, and $W^m(I)=(W_1(I),\dots, W_m(I), 0, 0,\dots)$, for any $m\in \mathbb{N}$.
       
       Fix  $n_0\in \mathbb{N}$. By (2.2),      for any $\rho>0$, there exists $m_0\in \mathbb{N}$ , depending only on $n_0$ and $\rho$, such that if $m\geqslant m_0$, then 
     \begin{equation}|F_k(I,\varphi)-F_k(I^{m},\varphi^{m})|\leqslant \rho,\quad \forall (I,\varphi)\in B^I_p(M_1)\times\mathbb{T}^{\infty},\label{dn}
     \end{equation}
     where $k=1,\cdots, n_0$.

             From now on, we always assume that 
             \[(I,\varphi)\in B_p^I(M_1)\times\mathbb{T}^{\infty},\quad\mbox{i.e.}\quad v\in B_p^v(\sqrt{M_1}).\] 
              By Lemma 1.4, we have 
             \begin{equation}
             \eqalign{|G_j(I,\varphi)|\leqslant \frac{C_0(j,M_1)}{\sqrt{I_j}},\\
     |\frac{\partial F_k}{\partial I_j}(I,\varphi)|\leqslant \frac{C_0(k,j, M_1)}{\sqrt{I_j}},\\
      |\frac{\partial F_k}{\partial \varphi_j}(I,\varphi)|\leqslant C_0(k,j,M_1).\label{clumsy1}
         }
         \end{equation}
                     From Lemma 1.3 and Lemma 2.1, we know that
                     \begin{equation}
                     \eqalign{
                      |W_j(I)-W_j(\bar{I})|\leqslant C_1(j,M_1)|I-\bar{I}|_1,\\
                            |\langle F_k\rangle(I)-\langle F_k \rangle(\bar{I})| \leqslant C_1(k,j,M_1)|I-\bar{I}|_1.\label{clumsy2}
                            }
                            \end{equation}
                            By (2.1) we get
                            \begin{equation}
                            |F_k(I^{m_0},\varphi^{m_0})-F_k(\bar{I}^{m_0},\bar{\varphi}^{m_0})|\leqslant C_2(k,m_0,M_1)|v^{m_0}-\bar{v}^{m_0}|,\label{clumsy3}
                            \end{equation}
                            where $|\cdot|$ is the maximum norm.
                            
                            We denote
                            \[C_{M_1}^{n_0,m_0}=m_0\cdot\max\{C_0,C_1,C_2:1\leqslant j\leqslant m_0,1\leqslant k\leqslant n_0\}.\]
                 Below we define a number of sets, depending on various parameters. All of them also depend on $m_0$ and $n_0$, but this dependence is not indicated.            For any $\delta>0$, and $T_0>0$, we define a subset $E(\delta,T_0)\subset B_p^I(M_1)$ as the  collection of all $I\in B_p^I(M_1)$ such that for every $\varphi\in \mathbb{T}^{\infty}$ and any $T\geqslant T_0$, we have 
      \begin{equation}
      \Big|\frac{1}{T}\int_0^{T}[F_k(I^{m_0},\varphi^{m_0}+W^{m_0}(I)t)-\langle F_k\rangle(I^{m_0})]dt\Big|\leqslant\delta,
     \label{average}
      \end{equation}
      for $k=1,\cdots,n_0$. 
      Let $\mathcal{S}_\epsilon^t$ be the flow generated by regular solutions of the system (1.5). We define two more groups of sets.
      \[S(t)=S(t,\epsilon,\delta, T_0,I,\varphi):=\{t_1\in[0,t]:\mathcal{S}_\epsilon^{t_1}(I,\varphi)\notin E(\delta,T_0)\times\mathbb{T}^{\infty}\}.\]
      \[N(\tilde{T})=N(\tilde{T},\epsilon,\delta,T_0):=\{(I,\varphi)\in \Gamma_0: \mbox{Mes}[S(\epsilon^{-1}, \epsilon,\delta, T_0,I,\varphi)]\leqslant \tilde{T}\}.\]
      Here and below $\mbox{Mes}[\cdot]$ stands  for the  Lebesgue measure in $\mathbb{R}$. 
      
      Clearly, $E(\delta,T_0)$ is a closed subset of $B_p^I(M_1)$ and $S(t,\delta, T_0, I, \varphi)$ is a open subset of $[0,t]$. The following result is the main lemma of this work:     
          \medskip

      \noindent{\bf Lemma 4.1.} For $k=1,\dots,n_0$, the $I_k$-component of any regular solution of (\ref{iapkdv}) with initial data in $N(\tilde{T},\epsilon,\delta,T_0)$ can be written as:\[I_k(t)=I_k(0)+\epsilon\int_0^t\langle F_k\rangle(I(s))ds +\Xi(t),\]
      where for any $\gamma\in(0,1)$ the function $|\Xi(t)|$ is bounded  on $[0,\frac{1}{\epsilon}]$ by
     \[
     \eqalign{ 4\epsilon C_{M_1}^{n_0,m_0}\Big\{\Big[2(\gamma+2T_0C_{M_1}\epsilon)^{1/2}\Big](T_0+\tilde{T}+\epsilon^{-1})\\
     +\Big[\frac{ T_0C_{M_1}\epsilon}{\gamma^{1/2}} +T_0C_{M_1}\epsilon
                +(\frac{T_0\epsilon}{2\gamma^{1/2}}+\frac{\epsilon C_{M_1}T_0^2}{3})\Big](T_0+\tilde{T}+\epsilon^{-1})\Big\}\\
                + 2\epsilon C_{M_1} \tilde{T}
                +2\rho+2\delta + 2\epsilon C_{M_1}(T_0+\tilde{T}).
             }
             \]
      
     \noindent{\it Proof}:         \quad For any $(I,\varphi)\in N(\tilde{T})$, we consider the corresponding set $S(t)$. It is composed of open intervals of total length less   than $\min\{\tilde{T},t\}$. Thus at most $[\tilde{T}/T_0]$ of them have length greater than or equal to $T_0$. We denote these long intervals by $(a_i, b_i)$, $1\leqslant i\leqslant d$, $d\leqslant \tilde{T}/T_0$ and denote by 
       $C(t)$  the complement of $\cup_{1\leqslant i\leqslant d}(a_i,b_i)$ in $[0,t]$. 
      
     By (\ref{dn}), we have \[\int_0^t F_k(I(s),\varphi(s))dt=\int_{C(t)}F_k(I^{m_0}(s),\varphi^{m_0}(s))ds + \xi_1(t),\]
      where $ |\xi_1(t)|\leqslant C_{M_1} \tilde{T} +\rho t$.\\
    
      The set $C(t)$ is composed of segments $[b_{i-1},a_i]$ (if necessary, we set $b_0=0$, and $a_{d+1}=t)$.
      We proceed by dividing each segment $[b_{i-1},a_i]$ into shorter segments by points $t_j^i$, where  $b_i=t_1^i < t_2^i< \cdots <t_{n_i}^i=a_i$.  The points $t_j^i$ lie outside the set $S(t)$ and $T_0\leqslant t_{j+1}^i-t_j^i\leqslant 2 T_0$ except for the terminal segment containing the end points $a_i$, which may be shorter than $T_0$.
      \smallskip
      
  This partition is constructed as follows:
      
              \begin{enumerate}
      \item[-----]   If $a_i-b_{i-1}\leqslant 2T_0$, then we keep the whole segment with no subdivisions. ($t_1^i=b_{i-1}$ ,$ t_2^i=a_i$).
      \item[-----]  If $a_i-b_{i-1}>2T_0$, we divide the segment in the following way:
      {\begin{enumerate}
          \item[a)]  If $b_{i-1}+2T_0$ does not belong to $S(t)$, we chose $t_2^i=b_{i-1}+2T_0$, and continue by subdividing  $[t_2^i, a_i]$;
          \item[b)] if $b_{i-1}+2T_0$ belongs to $S(t)$, then there are points in $[b_{i-1}+T_0,b_{i-1}+2T_0]$  which do not, by definition of $b_{i-1}$. We set $t_2^i$ equal  to one of these points and continue by subdividing $[t_2^i,a_i]$.
          \end{enumerate}}
                
      \end{enumerate}
      
      We will adopt the notation: $h_j^i=t_{j+1}^i-t_j^i$ and $s(i,j)=[t_j^i,t_{j+1}^i]$. So \[C(t)=\bigcup_{i=1}^{d}\bigcup_{j=1}^{n_i-1}s(i,j),\;T_0\leqslant h_j^i=|s(i,j)|\leqslant 2T_0,\;j\leqslant n_i-2.\]
      
   By its definition, $C(t)$ contains at most $[\tilde{T}/T_0]+1$ segments $[b_{i-1}, a_i]$, thus $C(t)$ contains at most  $[\tilde{T}/T_0]+1$  terminal subsegments of length less than $T_0$.  Since all other segments have length no less than $T_0$ and $t\leqslant \frac{1}{\epsilon}$,     the number of these segments is not greater than $[\epsilon T_0]^{-1}$.  So the total number of subsegments $s(i,j)$ is bounded by $1+[(\epsilon T_0)^{-1}]+[\tilde{T}/T_0]$.
     
      For  each segment $s(i,j)$ we define a subset $\Lambda(i,j)$ of $\{1,2,\cdots,m_0\}$ in the following way:
      \[l\in \Lambda(i,j) \quad \Longleftrightarrow\quad\exists t\in s(i,j),\quad I_l(t)<\gamma.\]
       If $l\in \Lambda$, then by (\ref{boundf}) we have
        \begin{equation}|I_l(t)|<2T_0C_{M_1}\epsilon +\gamma, \quad t\in s(i,j).\label{boundi}
        \end{equation}
     For $I=(I_1,I_2,\dots)$ and $\varphi=(\varphi_1,\varphi_2,\dots)$ we set
      \[\lambda_{i,j}(I)=\hat{I},\quad
      \lambda_{i,j}(\varphi)=\hat{\varphi},
      \]
      where $\hat{\varphi}=(\hat{\varphi_1},\hat{\varphi}_2,\dots)$ and $\hat{I}=(\hat{I}_1,\hat{I}_2,\dots)$ are defined by  the following relation:
      \[\mbox{If}\quad l\in \Lambda(i,j),\quad\mbox{then}\quad \hat{I}_l=0,\;\;\hat{\varphi}_l=0,\quad \mbox{else}\quad \hat{I}_l=I_l,\;\;\hat{\varphi}_l=\varphi_l.\]
      We also denote $\lambda_{i,j}(I,\varphi)=(\lambda_{i,j}(I),\lambda_{i,j}(\varphi))$ and when the segment $s(i,j)$ is clearly indicated, we write for short $\lambda_{i,j}(I,\varphi)=(\hat{I},\hat{\varphi})$.      
      
      Then on $s(i,j)$, using (\ref{clumsy3}) and (\ref{boundi}) we obtain 
      \begin{equation}
      \eqalign{
       \int_{s(i,j)}\Big|F_k\Big(I^{m_0}(s),\varphi^{m_0}(s)\Big)-F_k\Big(\lambda_{i,j}\big(I^{m_0}(s),\varphi^{m_0}(s)\big)\Big)\Big|ds\\
       \leqslant\int_{s(i,j)}C^{n_0,m_0}_{M_1}\Big|I^{m_0}(s)-\lambda_{i,j}\Big(I^{m_0}(s)\Big)\Big|^{1/2}ds\\
       \leqslant 2T_0C_{M_1}^{n_0,m_0}(\gamma +2T_0C_{M_1}\epsilon)^{1/2}. \label{p0}
       }
      \end{equation}
      
      In Proposition 1-5 below, $k=1,\dots, n_0$.
      \smallskip
      
      \noindent{\bf Proposition 1. }\[\int_{C(t)}F_k\Big(I^{m_0}(s),\varphi^{m_0}(s)\Big)ds =\sum_{i,j}\int_{s(i,j)}F_k\Big(I^{m_0}(t_j^i),\varphi^{m_0}(s)\Big)ds +\xi_2(t),\]
      where 
      \begin{equation}|\xi_2|\leqslant 4C^{n_0,m_0}_{M_1}\Big[(\gamma+2T_0C_{M_1}\epsilon)^{1/2}+\gamma^{-1/2}T_0C_{M_1}\epsilon\Big](T_0+\tilde{T}+\epsilon^{-1}).\label{p1}
      \end{equation}
      
     \noindent{\it Proof}: We may write $\xi_2(t)$ as 
    \[\eqalign{
    \xi_2(t) =\sum_{i,j}\int_{s(i,j)}\Big[F_k\Big(I^{m_0}(s),\varphi^{m_0}(s)\Big)-F_k\Big(I^{m_0}(t_j^i),\varphi^{m_0}(s)\Big)\Big]ds\\
   \qquad :=\sum_{i,j} I(i,j).
    }\]
     For each $s(i,j)$, we have
     \begin{equation}
     \eqalign{\int_{s(i,j)}\Big|F_k\Big(\hat{I}^{m_0}(s),\hat{\varphi}^{m_0}(s)\Big)-F_k\Big(\hat{I}^{m_0}(t_j^i),\hat{\varphi}^{m_0}(s)\Big)\Big|ds\\
     \leqslant \int_{s(i,j)}\gamma^{-1/2}C^{n_0,m_0}_{M_1}\Big|\hat{I}^{m_0}(s)-\hat{I}^{m_0}(t_j^i)\Big|ds\\
     \leqslant 2\gamma^{-1/2}T_0^2C_{M_1}\epsilon.
     }
     \label{p1.1}
     \end{equation}
        We replace the integrand $F_k(I^{m_0},\varphi^{m_0})$ by $F_k(\hat{I}^{m_0},\hat{\varphi}^{m_0})$. Using (\ref{p0})  and (\ref{p1.1}) we obtain that
        \[I(i,j)\leqslant 4T_0C^{n_0,m_0}_{M_1}\Big[(\gamma+2T_0C_{M_1}\epsilon)^{1/2}+\gamma^{-1/2}T_0C_{M_1}\epsilon\Big].\]        The inequality (\ref{p1}) follows.\quad $\square$
        \medskip

    On each subsegment $s(i,j)$, we now consider  the unperturbed linear dynamics  $\varphi_j^i(t) $ of the angles $\varphi^{m_0}\in \mathbb{T}^{m_0}$ :
    \[\varphi_j^i(t)=\varphi^{m_0}(t_j^i)+W^{m_0}(I(t_j^i))(t-t_j^i)\in\mathbb{T}^{m_0},\quad t\in s(i,j).\]
    \smallskip
    
    \noindent{\bf Proposition 2.}\[\sum_{i,j}\int_{s(i,j)}F_k\Big(I^{m_0}(t_j^i),\varphi^{m_0}(s)\Big)ds=\sum_{i,j}\int_{s(i,j)}F_k\Big(I^{m_0}(t_j^i),\varphi_j^i(s)\Big)ds +\xi_3(t),\]
    where
     \begin{equation}
     \eqalign{|\xi_3(t)|\leqslant 4C_{M_1}^{n_0,m_0}(\gamma+2T_0C_{M_1}\epsilon)^{1/2}(T_0+\tilde{T}+\epsilon^{-1})\\
    +(C^{n_0,m_0}_{M_1})^2\Big( \frac{2T_0\epsilon}{\gamma}+\frac{4\epsilon C_{M_1}T_0^2}{3}\Big)(T_0+\tilde{T}+\epsilon^{-1}).
    }
    \label{p2}
    \end{equation}
 
 \noindent{\it Proof}:  For  each $s(i,j)$ we have 
 \[
 \eqalign{\int_{s(i,j)}\Big|\lambda_{i,j}\Big(\varphi^{m_0}(s)-\varphi_j^i(s)\Big)\Big|ds\\
       \leqslant \int_{s(i,j)}\int_{t_j^i}^{s}\Big| \lambda_{i,j}\Big(\epsilon G^{m_0}(I(s^{\prime}),\varphi(s^{\prime}))+W^{m_0}(I(s^{\prime}))-W^{m_0}(I(t_j^i))\Big )\Big|ds^{\prime}ds\\
    \leqslant \int_{s(i,j)}\int_{t_j^i}^{s}C_{M_1}^{n_0,m_0}\Big[\epsilon \gamma^{-1/2}+|I(s^{\prime})-I(t_j^i))|_1\Big]ds^{\prime} ds\\
    \leqslant \int_{s(i,j)}C_{M_1}^{n_0,m_0}\Big[\gamma^{-1/2}\epsilon(s-t_j^i)+\frac{1}{2}C_{M_1}\epsilon(s-t_j^i)^2\Big]ds\\
    \leqslant C_{M_1}^{n_0,m_0}\Big(\frac{2T_0^2\epsilon}{\sqrt{\gamma}} +\frac{4\epsilon C_{M_1} T_0^3}{3}\Big).
    }
    \]
    Here the first inequality comes from equation (\ref{angle1}),  and using (\ref{clumsy1}) and (\ref{clumsy2}) we can get the second inequality. The third one follows from (\ref{boundf}).
    
    Using again (\ref{clumsy1}), we get
 \[\eqalign{\int_{s(i,j)}\Big[F_k\Big(\lambda_{i,j}\big(I^{m_0}(t_j^i),\varphi^{m_0}(s)\big)\Big)-F\Big(\lambda_{i,j}\big(I^{m_0}(t_j^i),\varphi_j^i(s)\big)\Big)\Big]ds\\   
       \leqslant \int_{s(i,j)}C_{M_1}^{n_0,m_0}\Big|\lambda_{i,j}\Big(\varphi^{m_0}(s)-\varphi_j^i(s)\Big)\Big|ds\\
    \leqslant (C_{M_1}^{n_0,m_0})^2\Big(\frac{2T_0^2\epsilon}{\sqrt{\gamma}} +\frac{4\epsilon C_{M_1} T_0^3}{3}\Big).
    }
    \]
    Therefore (\ref{p2}) holds for the same reason as (\ref{p1}).\quad\quad$\square$
    \smallskip
      
      We will now compare the integral $\int_{s(i,j)}F_k(I^{m_0}(t_j^i),\varphi_j^i(s))ds$ with the average value $\langle F_k(I^{m_0}(t_j^i))\rangle h_j^i$.
      \smallskip

     \noindent{\bf Proposition 3.}\[\sum_{i,j}\int_{s(i,j)}F_k\Big(I^{m_0}(t_j^i),\varphi_j^i(s)\Big)ds=\sum_{i,j}h_j^i\langle F_k\rangle\Big(I^{m_0}(t_j^i)\Big)+\xi_4(t),\]
      where
      \begin{equation}|\xi_4(t)|\leqslant \frac{2\delta}{\epsilon}+2C_{M_1}(T_0+\tilde{T}).\label{p3}
      \end{equation}
      
      \noindent{\it Proof}: \quad
      We divide the set of segments $s(i,j) $ into two subsets $\Delta_1 $ and $\Delta_2$. Namely, $s(i,j)\in \Delta_1$  if $h^i_j\geqslant T_0$ and $s(i,j)\in \Delta_2$ otherwise.
      
      (i) $s(i,j)\in \Delta_1$. In this case, by (\ref{average}), we have
      \[\Big|\int_{s(i,j)}\Big[F_k\Big(I^{m_0}(t_j^i),\varphi_j^i(s)\Big)-\langle F_k\rangle\Big(I^{m_0}(t_j^i)\Big)\Big]ds\Big|\leqslant \delta h_j^i.\]
      So\[\sum_{s(i,j)\in\Delta_1}\Big|\int_{s(i,j)}F_k\Big(I^{m_0}(t_j^i),\varphi_j^i(s)\Big)ds-\langle F_k\rangle\Big(I^{m_0}(t_j^i)\Big)h_j^i\Big|\leqslant \delta \sum_{s(i,j)\in\Delta_1} h_j^i \leqslant \frac{2\delta}{\epsilon}.
      \]

      (ii)  $s(i,j)\in\Delta_2$. Now, using (\ref{boundf}) we get
          \[\Big|\int_{s(i,j)}F_k\Big(I^{m_0}(t_j^i),\varphi_j^i(s)\Big)ds - \langle F_k\rangle\Big(I^{m_0}(t_j^i)\Big)h_j^i\Big|\leqslant 2C_{M_1} h_j^i\leqslant 2C_{M_1}T_0.
          \]
      Since $\mbox{Card}(\Delta_2)\leqslant (1+\tilde{T}/T_0)$, then
      
      \[\sum_{s(i,j)\in\Delta_2}\Big|\int_{s(i,j)}F\Big(I^{m_0}(t_j^i),\varphi_j^i(s)\Big)ds - \langle F_k\rangle\Big(I^{m_0}(t_j^i)\Big)h_j^i\Big|\leqslant 2C_{M_1}(\tilde{T}+T_0).
      \]
      This implies the inequality (\ref{p3}).\quad $\square$
      \smallskip

      \noindent{\bf Proposition 4.}
      \[\sum_{i,j}h_j^i\langle F_k\rangle\Big(I^{m_0}(t_j^i)\Big)=\int_{C(t)}\langle F_k\rangle\Big(I^{m_0}(s)\Big)ds +\xi_5(t),\]
      where\begin{equation}
      |\xi_5(t)|\leqslant4\epsilon C_{M_1}C_{M_1}^{n_0,m_0}T_0(T_0+\tilde{T}+\epsilon^{-1}).\label{p4}
      \end{equation}  
         
     \noindent{\it Proof}: \quad  Indeed, as \[|\xi_5(t)|=\Big|\sum_{i,j}\int_{s(i,j)}\Big[\langle F_k\rangle\big(I^{m_0}(s)\big)-\langle F_k\rangle\big(I^{m_0}(t_j^i)\big)\Big]ds\Big|,\]
      
     using  (\ref{boundf}) and (\ref{clumsy2}) we get
     \[
     \eqalign{|\xi_5(t)|&\leqslant \sum_{i,j}\int_{s(i,j)}C_{M_1}^{n_0,m_0}|I^{m_0}(s)-I^{m_0}(t^i_j)|ds \\
     &\leqslant\epsilon \sum_{i,j}C_{M_1}C_{M_1}^{n_0,m_0}(h_j^i)^2\leqslant 4\epsilon C_{M_1}C_{M_1}^{n_0,m_0}T_0(T_0+\tilde{T}+\epsilon^{-1}).\quad \square 
     }\]
    
     Finally,

       \noindent{\bf Proposition 5.} \[\int_{C(t)}\langle F_k\rangle\Big(I^{m_0}(s)\Big)ds=\int_0^t\langle F_k\rangle\Big(I(s)\Big)ds +\xi_6(t),\]  
      and $|\xi_6(t)|$ is bounded by $C_{M_1}\tilde{T}+\rho t$.\quad $\square$
      \medskip

      Gathering the estimates in Propositions 1-5, we obtain
       \[
       \eqalign{ I_k(t)&=I_k(0)+\epsilon \int_0^t F_k\Big(I(s),\varphi(s)\Big)ds
       \\&=I_k(0)+ \epsilon \int_0^t \langle F_k\rangle\Big(I(s)\Big)ds +\Xi(t),
       }\]
      where \[
      \eqalign{|\Xi(t)|&\leqslant\epsilon \sum_{i=1}^6|\xi_i(t)|\\
                & \leqslant 4\epsilon C_{M_1}^{n_0,m_0}\Big[2(\gamma+2T_0C_{M_1}\epsilon)^{1/2}+\frac{ T_0C_{M_1}\epsilon}{\gamma^{1/2}} +T_0C_{M_1}\epsilon\\
                &\quad+\Big(\frac{T_0\epsilon}{2\gamma^{1/2}}+\frac{\epsilon C_{M_1}T_0^2}{3}\Big)\Big](T_0+\tilde{T}+\epsilon^{-1})
                + 2\epsilon C_{M_1} T_1\\
                &\quad+2\rho+2\delta + 2\epsilon C_{M_1}(T_0+\tilde{T}),
                \quad t\in[0, \frac{1}{\epsilon}].
             }\]

      Lemma 4.1 is proved.\quad $\square$\\
      
      \noindent{\bf Corollary 4.2.} For any $\bar{\rho}>0$, with  a suitable choice of  $\rho$, $\gamma$, $\delta$, $T_0$, $\tilde{T}$, the function $|\Xi(t)|$ in Lemma 4.1 can be made smaller than $\bar{\rho}$, if $\epsilon$ is small enough.
      
      \noindent{\it Proof}:\quad
      We choose
      \[\gamma=\epsilon^{\alpha},\; T_0=\epsilon^{-\sigma},\; \tilde{T}=\frac{\bar{\rho}}{9C_{M_1}\epsilon},\;\delta=\rho=\frac{\bar{\rho}}{9}\]
      with\[1-\frac{\alpha}{2}-\sigma>0,\;0<\sigma<\frac{1}{2}.\]Then for $\epsilon$ sufficiently small we have \[|\Xi(t)|<\bar{\rho}.\quad \square\]
      \smallskip

        On the Hilbert space $h^p$, we adopt a $\zeta_0$-admissible Gaussian measure $\mu$.       Define  corresponding measures $\mu_I=\Pi_I\circ\mu$ and $\mu_{I,\varphi}=\Pi_{I,\varphi}\circ\mu$ in the  spaces $h_{I+}^p$ and  $h_{I+}^p\times \mathbb{T}^{\infty}$.

      \noindent{\bf Lemma 4.3.} The measure $\mu_{I,\varphi} $ is  a product measure  $d\mu_{I,\varphi}=d\mu_I d\varphi$, where $d\varphi$ is the Haar measure on $\mathbb{T}^{\infty}$.
      
      \noindent {\it Proof}: Since the measure $\mu$ is invariant under rotations $\Phi_{\theta}$, the $\Pi_{\varphi}\circ d\mu$ is a measure on $\mathbb{T}^{\infty}$, invariant under the rotations. So this is the Haar measure $d\varphi$. Consequently  the image of the measure $\mu_{I,\varphi}$ under the natural projection $(I,\varphi)\mapsto~\varphi$ is $d\varphi$. Since the spaces $h^p_{I+}$ and $\mathbb{T}^{\infty}$ are separable, then for $\varphi\in~\mathbb{T}^{\infty}$  there exists a Borel probability measure $\pi_{\varphi}(dI)$ on $h^p_{I+}$ such that $\mu_{I,\varphi}=~\pi_{\varphi}(dI)d\varphi$. That is, for any bounded continuous function $f(I,\varphi)$, we have 
      \[\langle \mu_{I,\varphi}, f\rangle=\int_{\mathbb{T}^{\infty}}\Big(\int_{h^p_{I+}}f(I,\varphi)\pi_{\varphi}(dI)\Big)d\varphi. 
        \]
        (see e.g.  \cite{13}). For any $\theta\in\mathbb{T}^{\infty}$ we have 
        \[
        \eqalign{
        \langle \mu_{I,\varphi},f \rangle&=\langle \mu_{I,\varphi},f\circ\Phi_{\theta}\rangle \\
        &=\int\int  f(I,\varphi+\theta)\pi_{\varphi}(dI)d\varphi=\int\int f(I,\varphi)\pi_{\varphi-\theta}(dI)d\varphi.
        }
        \]
        Integrating in $d\theta$ we see that 
        \[\mu_{I,\varphi}(dId\varphi)=d\mu^{\prime}(dI)d\varphi,\]
        where $d\mu^{\prime}(dI)=\int_{\mathbb{T}^{\infty}}\pi_{\theta}(dI)d\varphi$. We must have $d\mu^{\prime}=d\mu_I$, and the assertion of the lemma is proved.\quad $\square$
        \smallskip
        
        The two lemmas below deal with the sets $E$ and $N$, defined at the beginning of this section.
        \smallskip
        
       \noindent {\bf Lemma 4.4.} For any $\delta>0$, $\lim_{T_0\rightarrow \infty}\mu_I(B^I_p(M_1)\setminus E(\delta,T_0))=0$.
       
        \noindent{\it Proof}: \quad From the definition of $E(\delta, T_0)$, we know that  \[E(\delta, T_0)\subset E(\delta,T_0^{\prime}),\quad \mbox{if}\quad T_0\leqslant T_0^{\prime}.\]
        Let 
        $E_{\infty}(\delta):=\bigcup_{T_0>0}E(\delta,T_0).$
        Due to the inclusion above we have to check that 
        \[\mu_I(B^I_p(M_1)\setminus E_{\infty}(\delta))=0.\]
        Denote\[\mathcal{R}(N):=\bigcup_{L\in \mathbb{Z}^{m_0}\setminus\{0\},|L|\leqslant N}\{I\in B^I_p(M_1):\quad W^{m_0}(I)\cdot L=0\},
        \]
        where $W^{m_0}(I)=(W_1(I),\dots,W_{m_0}(I))$.
        Let us write $F_k(I^{m_0},\varphi^{m_0})$ as a  Fourier series $F_k(I^{m_0},\varphi^{m_0})=\sum_{L\in\mathbb{Z}^{m_0}} F_k^Le^{iL\cdot \varphi^{m_0}}$, where $F_k^L=F_k^L(I^{m_0})$.
        Then there exists $N_0>0$ such that 
       \[ \Big|F_k(I^{m_0},\varphi^{m_0})-\sum_{|L|\leqslant N_0} F_k^Le^{iL\cdot \varphi^{m_0}}\Big|<\frac{\delta}{2},\quad  k=1,\cdots, n_0.\]
        Arguing as in the proof  of  Lemma 2.2, we see that if $I\notin \mathcal{R}(N_0)$,  then
         \[\Big|\sum_{0\neq|L|\leqslant  N_0} \frac{1}{T_0}\int_0^{T_0} F_k^Le^{iL\cdot W^{m_0}t}dt\Big|\leqslant \frac{2}{T_0}\Big(\inf_{0\neq|L|\leqslant N_0}|L\cdot W^{m_0}|\Big)^{-1}\sum_{|L|\leqslant N_0}|F_k^L|.
         \]
         where $W^{m_0}=W^{m_0}(I)$. 
         The r.h.s of the above inequality can be made smaller  than $\delta/2$ by choosing  $T_0$  large enough. So 
         we have \[B^I_p(M_1)\setminus\mathcal{R}(N_0)\subset E_{\infty}(\delta),\]
         and it remains to show that 
         \[\mu_I(\mathcal{R}(N_0))=0.\]
     By  Lemma 1.2, \[W^{m_0}(I)\cdot L \not\equiv 0,\quad\forall L\in \mathbb{Z}^{m_0}\setminus\{0\},\]
        Since $W^{m_0}(I)$ is analytic with respect to $I$  and $\mu_I$ is a non-degenerated Gaussian measure, then due to  Theorem 1.6 in \cite{8}, for any $L\in \mathbb{Z}^{m_0}$,  we have\[\mu_I(\{I\in h^p_I: W^{m_0}(I)\cdot L=0\})=0.\]
        Therefore,
        \[\mu_I(\mathcal{R}(N_0))=0.\quad \square\]        
              \smallskip

        \noindent{\bf Lemma 4.5.} Fix any $\delta>0$, $\bar{\rho}>0$.  Then for every $\nu>0$ we can find $T_0>0$ such that \[\mu_{I,\varphi}\Big(\Gamma_0\setminus N\Big)< \nu,\] where $N=N(\frac{\bar{\rho}}{9C_{M_1}\epsilon},\epsilon,\delta, T_0)$.

      \noindent{\it Proof}:\quad Let us denote $\Gamma_E=E(\delta,T_0)\times\mathbb{T}^{\infty}$, $\Gamma_1=B^I_p(M_1)\times\mathbb{T}^{\infty}$ and $\Gamma_E^{\infty}:=\bigcup_{T_0>0}\Gamma_E(\delta,T_0)$.
      Since the sets $\Gamma_E(\delta,T_0)$ are increasing with $T_0$, then 
 from Lemmas 4.3 and 4.4 we know that 
 \begin{equation}
 \lim_{T_0\to\infty}\mu_{I,\varphi}(\Gamma_1\setminus\Gamma_E(\delta,T_0))=\mu_{I,\varphi}(\Gamma_1\setminus\Gamma_E^{\infty})=0.\label{4.16}
 \end{equation}
       Let  $d\mu_1$ be the measure $d\mu dt$ on $h^p\times \mathbb{R}$, and $\mathcal{S}_{v,\epsilon}^t$ be the flow of the perturbed KdV equation (\ref{tvpkdv}) on $h^p$.
      We now define following  subset of $h^p\times \mathbb{R}$:
      \[B^{\prime}=\Big\{(v,t):\mathcal{ S}_{v,\epsilon}^t( v)\in \Pi_{I,\varphi}^{-1}(\Gamma_1\setminus\Gamma_E(\delta,T_0)),v\in B_p^v(\sqrt{M_0}),t\in[0,\frac{1}{\epsilon}]\Big\}.\] 
   By Theorem 3.6, there exists a constant $C_2(M_1)$ depending only on  $M_1$ such that
      \[
      \eqalign{\mu_1(B^{\prime})=& \int_0^{\epsilon^{-1}}\mu\Big(\mathcal{S}_{v,\epsilon}^{-t}\Big(\Pi_{I,\varphi}^{-1}(\Gamma_1\setminus\Gamma_E(\delta,T_0)\Big)\bigcap \Pi_{I,\varphi}^{-1}(\Gamma_0)\Big)dt\\
      &\leqslant \frac{1}{\epsilon} e^{C_2(M_1)} \mu\Big(\Pi_{I,\varphi}^{-1}(\Gamma_1\setminus\Gamma_E(\delta,T_0))\Big)\\
      &=\frac{1}{\epsilon} e^{C_2(M_1)} \mu_{I,\varphi}\Big(\Gamma_1\setminus\Gamma_E(\delta,T_0)\Big).
           }\label{mu1}
      \]
       For $v\in \Pi_{I,\varphi}^{-1}(\Gamma_0)$, we define \[S(I,\varphi)=S(v)=\{ t\in[0,\epsilon^{-1}]: \mathcal{S}_{v,\epsilon}^t(v)\in B_p^v(\sqrt{M_1})\setminus\Pi_{I,\varphi}^{-1}(\Gamma_E(\delta, T_0))\}.\] 
      By the Fubini theorem, we have \[\mu_1(B^{\prime})=\int_{\Pi_I^{-1}(\Gamma_0) }\mbox{Mes}(S(v))\mu(dv),\]
     Thus
   \[
   \eqalign{\mu_{I,\varphi}(\Gamma_0\setminus N)=\mu_{I,\varphi}\Big(\{(I,\varphi)\in\Gamma_0:\mbox{Mes}(S(I,\varphi)>\frac{\bar{\rho}}{9C_{M_1}\epsilon}\}\Big)\\
   \leqslant \frac{9C_{M_1}e^{C_2(M_1)}}{\bar{\rho}}\mu_{I,\varphi}\Big(\Gamma_1\setminus\Gamma_E(\delta,T_0)\Big),
   }
   \]
   by the Chebyshev inequality. In view of (\ref{4.16})             the term on  the right hand side becomes arbitrary small when $T_0$ is large enough.
   The statement of Lemma 4.5 follows.\quad $\square$
   \smallskip

  We pass to   the slow time $\tau=\epsilon t$. Let  $v^{\epsilon}(\tau)$, $\tau\in[0,1]$,  be a solution of the equation (\ref{vpkdv}) and $(I^{\epsilon}(\tau),\varphi^{\epsilon}(\tau))=\Pi_{I,\varphi}(v^{\epsilon}(\tau))$.
   
    By Lemma 2.1 and (\ref{prop 3.1}), we know that for any $p\geqslant 0$,  the mapping      
            \[F_J: h^p_I\to h^{p+\zeta_0}_I,\quad J\mapsto\langle F\rangle(J),\]
            where $\langle F\rangle(J)=(\langle F_1\rangle(J),\langle F_2\rangle(J),\dots)$  is analytic. Hence, there exists  $C_3(M_1)$ such that
            \begin{equation}|F_J(J_1)-F_J(J_2)|_p\leqslant C_3(M_1)|J_1-J_2|_p,\quad J_1,J_2\in B_p^I(2M_1).\label{l1}
            \end{equation}
            Using Picard's theorem, for any $J_0\in B^I_p(M_1)$ there exists  a unique solution $J(t)$ of the averaged equation (\ref{apkdv1}) with $J(0)=J_0$. 
                We denote \[T(J_0):=\inf\{\tau>0: |J(\tau)|_p>2M_1\}.\]

          Now we are in a position to prove the assertion (i) of Theorem 0.2.
          \smallskip
                   
                   For any $\bar{\rho}>0$,  there exist $n_1$ such that 
            \begin{equation}
            \eqalign{|F(I,\varphi)-F^{n_1}(I,\varphi)|_p<\frac{\bar{\rho}}{8}e^{-C_3(M_1)},\quad (I,\varphi)\in B^I_p(2M_1)\times\mathbb{T}^{\infty},\\
            |\langle F\rangle(J)-\langle F\rangle^{n_1}(J)|_p<\frac{\bar{\rho}}{8}e^{-C_3(M_1)},\quad J\in B^I_p(2M_1).\label{l2}
            }
            \end{equation}
            Choose $\rho_0$ such that 
           
            \[8\sum_{j=1}^{n_1}j^{1+2p}\rho_0=\bar{\rho}e^{-C_3(M_1)}.\]
            
            By Lemmata 4.1 and 4.2, there is a set $\Gamma_{\bar{\rho}}=N(\frac{\rho_0}{9C_{M_1}\epsilon},\epsilon,\frac{\rho_0}{9},\epsilon^{-\sigma})$, $\sigma<1/2$,  such that if $\epsilon$ is small enough and  $(I^{\epsilon}(0),\varphi^{\epsilon}(0))\in \Gamma_{\bar{\rho}}$, then 
            \[
            I_k(\tau)=I_k(0)+\int_0^\tau \langle F_k\rangle (I(s))ds+\xi_k(\tau),\quad |\xi_k(\tau)|<\rho_0,\quad \tau\in[0,1],
                       \]
                       for $k=1,\cdots, n_1$.
            Therefore, by (\ref{l1}) and (\ref{l2}),
            \[|I(\tau)-J(\tau)|_p\leqslant \int_0^{\tau} C_3(M_1)|I(\tau)-J(\tau)|_p ds+\xi_0(\tau),\quad |\xi_0(\tau)|\leqslant \frac{\bar{\rho}}{2}e^{-C_3(M_1)},\]
            for $(I(0),\varphi(0))\in \Gamma_{\bar{\rho}}$, $I(0)=J(0)$ and $|\tau|\leqslant \min\{1,T(J(0))\}$. By Gronwall's lemma, 
           \[ |I(\tau)-J(\tau)|_p\leqslant \bar{\rho},\quad |\tau|\leqslant \min\{1,T(J(0))\}.\]
           Assuming  that $\bar{\rho}<<M_1$, we get from the definition of $T(J(0))$ that $T(J(0))$  is bigger than 1.
           This establishes inequality (\ref{main1}). From Lemma 4.5  we know that $\lim_{\epsilon\to0}\mu_{I,\varphi}(\Gamma_0-\Gamma_{\bar{\rho}})=0$.\quad $\square$
            
 \subsection{Proof of the assertion (ii)}

            It is not hard to see that the assertion for any $0\leqslant \bar{T}_1<\bar{T}_2\leqslant 1$ would follow if we can prove it for $\bar{T}_1=0$, $\bar{T}_2=1$. So we assume that $\bar{T}_1=0$, and $\bar{T}_2=1$. For any $(m,n)\in \mathbb{N}^2$, we fix $\alpha<1/8$, and denote 
            \[
            \eqalign{\mathcal{B}_m(\epsilon):=\Big\{ I\in B^I_p(M_1): \inf_{k\leqslant m}|I_k|<\epsilon^{\alpha}\Big\},\\
            \mathcal{R}_{m,n}(\epsilon):=\bigcup_{|L|\leqslant n,L\in\mathbb{Z}^m\setminus\{0\}}\Big\{I\in B^I_p(M_1):| W(I)\cdot L|<\epsilon^{\alpha}\}\Big\}.
            }\]
            Then let
            \[\Upsilon_{m,n}(\epsilon)=\Big(\bigcup_{m_0\leqslant m}\mathcal{R}_{m_0,n}(\epsilon)\Big)\cup\mathcal{B}_m(\epsilon).\]
           Denote 
            \[S(\epsilon, m, n, I_0,\varphi_0)=\{\tau\in[0,1]: I^{\epsilon}(\tau)\in\Upsilon_{m,n}(\epsilon)\}，\]
         and fix any $\nu>0$. Then using Theorem 3.6 and  arguing as in Lemma 4.4 and Lemma 4.5, we get that, for any $(m,n)\in \mathbb{N}^2$, there exists open subset $\Gamma_{\nu}^{m,n}\subset \Gamma_0$, $\epsilon_{m,n}>0$ and a positive function $\rho_{m,n}(\epsilon)$, converging to zero as~$\epsilon\rightarrow 0$,  such that 
            \[\mu_{I,\varphi}(\Gamma_0-\Gamma_{\nu}^{m,n})<\frac{\nu}{2^{mn}}\quad \mbox{and}\quad
            \mbox{Mes}(S(\epsilon,m,n,I_0,\varphi_0))\leqslant\rho_{m,n}(\epsilon), \]
            if $(I_0,\varphi_0)\in\Gamma_{\nu}^{m,n}$ and $\epsilon\leqslant\epsilon_{m,n}$.
            Let
            \[\Gamma_{\nu}=\bigcap_{(m,n)\in\mathbb{N}^2}\Gamma_{\nu}^{m,n},\]
            then
            \begin{equation}\mu_{I,\varphi}(\Gamma_0-\Gamma_{\nu})<\nu.\label{i}
            \end{equation}
            The sets $\Gamma_{\nu}$ may be chosen in such a manner  that
            \begin{equation}
            \Gamma_{\nu_1}\subset\Gamma_{\nu_2},\quad \mbox{if} \quad \nu_2<\nu_1.\label{ii}
            \end{equation}
           For any $(I_0,\varphi_0)\in \Gamma_{\nu}$,  consider a solution $(I^{\epsilon}(\tau),\varphi^{\epsilon}(\tau))$ such that \[(I^{\epsilon}(0),\varphi^{\epsilon}(0))=(I_0,\varphi_0).\]
           Fix  $m\in\mathbb{N}$, take a bounded Lipschitz function $g$ defined on the torus $\mathbb{T}^m\subset \mathbb{T}^{\infty}$ such that $Lip(g)\leqslant 1$ and $|g|_{L_{\infty}}\leqslant 1$. Let
           $\sum_{s\in\mathbb{Z}^{m}}g_s e^{i s\cdot \varphi}$ be  its Fourier series.  Then for any $\rho>0$, there exists $n$, such that if we denote $\bar{g}_n=\sum_{|s|\leqslant n}g_s e^{i s\cdot\varphi}$, then
           \[\Big|g(\varphi)-\bar{g}_n(\varphi)\Big|<\frac{\rho}{2},\quad\forall \varphi\in\mathbb{T}^{m}.\]
           
                 \quad For any $(I_0,\varphi_0)\in \Gamma_{\nu}$, we consider the set $S(\epsilon,m,n, I_0,\varphi_0)$. It is composed of open intervals of total length less   than $\tilde{T}=\rho_{m,n}(\epsilon)$.        
      Proceeding  as in Lemma 4.1 and Corollary 4.2, we find that
for $\epsilon$ small enough we have 
\[\Big|\int_0^1g\big(\varphi^{\epsilon, m}(\tau)\big)d\tau-\int_{\mathbb{T}^m}g(\varphi)d\varphi\Big|<\rho.\]

That is ,
\begin{equation}
\Big|\int g(\varphi)\mu_{\bar{T}_1,\bar{T}_2}^{\epsilon}(d\varphi)-\int g(\varphi) d\varphi\Big|\longrightarrow0\quad\mbox{as}\quad\epsilon\rightarrow 0,
\end{equation}
for any Lipschitz function as above. Hence, $\mu_{\bar{T}_1,\bar{T}_2}^{\epsilon}$ converges weakly to $d\varphi$ (see \cite{13}). This proves the required assertion with $\Gamma_\varphi$ replaced by $\Gamma_{\nu}$. Let us choose 
\[\Gamma_{\varphi}=\bigcup_{\nu>0}\Gamma_{\nu}.\]

Then 
\[\mu_{I,\varphi}(\Gamma_0-\Gamma_{\varphi})=0,\]
by (\ref{i}) and (\ref{ii}), and for any $(I_0,\varphi_0)\in \Gamma_{\varphi}$ the required convergence of measures holds. 
This proves the second assertion of Theorem 0.2.\quad $\square$

  \section{Application to a special case }
      
      In this section we prove Proposition 0.3.
       Clearly, we only need to prove the statement~(ii) of  assumption A. 
      Let $\mathcal{F}: H^m\rightarrow \mathbb{R}$ be a smooth functional (for some $m\geqslant 0$). If $u(t)$ is a solution of (\ref{pkdv}), then\[\frac{d}{dt}\mathcal{F}(u(t))=\langle\nabla \mathcal{F}(u(t)), -V(u)+\epsilon f(x)\rangle.\]
     In particular, if $\mathcal{F}(u)$ is an integral of motion for the KdV equation, then we have $\langle\nabla \mathcal{F}(u(t), V(u)\rangle=0$, so\[\frac{d}{dt}\mathcal{F}(u(t))=\epsilon\langle\nabla \mathcal{F}(u(t)), f(x)\rangle.\]
   Since $ ||u(0)||_0^2$  is an integral of motion, then 
   \[\frac{d}{dt}||u(t)||_0^2=2\epsilon\langle u,f(x)\rangle\leqslant\epsilon(||u||_0^2+||f(x)||_0^2).  \]  
        Thus we have \begin{equation}||u(t)||_0^2\leqslant e^{\epsilon t}(||u(0)||_0^2+\epsilon t||f(x)||_0^2).\label{(A.1)}\end{equation}

       The KdV equation has infinitively many integral of motion $\mathcal{J}_m(u)$, $m\geqslant 0$. The integral $\mathcal{J}_m$  can be writen as 
       \[\mathcal{J}_m(u)=||u||_m^2 +\sum_{r=3}^m \sum_{\mathbf{m}}\int C_{r,\mathbf{m}}u^{(m_1)} \cdots u^{(m_r)} dx,\]
        where the inner sum is taken over all integer $r$-vectors $\mathbf{m}=(m_1,\dots, m_r)$, such that $0\leqslant m_j \leqslant m-1$, $j=1,\dots, r$ and $m_1+\cdots+m_r= 4 + 2m - 2r$. Particularly, $\mathcal{J}_0(u)=||u||_0^2$.
        
        Lets consider\[I=\int u^{(m_1)}\cdots f^{(m_i)}\cdots u^{(m_{r_1})} dx, \quad m_1+\cdots+m_{r_1}=M,\]
        where $r_1\geqslant 2$, $M\geqslant 1$, and $0\leqslant m_j \leqslant \mu-1$.
        Then, by H\"older's inequality,
        \[|I|\leqslant ||u^{(m_1)}||_{L_{p_1}} \cdots||f(x)||_{L_{p_i}}\cdots ||u^{m_{r_1}}||_{L_{p_f}}, \quad p_j=\frac{M}{m_j}\leqslant \infty. \]
        Applying next the Gagliardo-Nirenberg  and the Young inequalities,  we obtain that \begin{equation}|I|\leqslant \delta ||u||_{\mu}^2 +C_{\delta}||u||_0^{C_1}, \quad \forall \delta >0,\label{(A.2)}
        \end{equation}
       where $C_{\delta}$ and $ C_1$ do not depend on $u$. Below we denote $C$ a positive constant independent of $u$,  not necessary the same in each inequality.
       Let\[I_1:=\langle\nabla \mathcal{J}_m(u), f\rangle=\langle u^{(m)},f^{(m)}\rangle+\sum_{r=3}^{m}\sum_{\mathbf{m}}C_{r,\mathbf{m}}^{'} u^{(m_1)}\cdots f^{(m_i)}\cdots u^{m_r} dx, \]
        where $m_1+\cdots +m_r=6+2m-2r$. Using (\ref{(A.2)}) with a suitable $\delta$, we get
        \begin{equation}I_1\leqslant ||u||_{m}^2 + C||u||_0^{C_1}\leqslant ||u||_m^2 +C(1+||u||_0^{4m})+||f||_m^2. \end{equation}

        If $u(t)=u(t,x)$ is a solution of equation (\ref{pkdv}), then
        \[\frac{d}{dt}\mathcal{J}_m(u)=\langle\nabla \mathcal{J}_m(u), \epsilon f\rangle
         \leqslant \epsilon ||u||_{m}^2+\epsilon C(1+||u||_0^{4m})+\epsilon ||f||_m^2,
         \]
    and     \[\frac{1}{2}||u||_m^2-C(1+||u||_0^{4m})\leqslant \mathcal{J}_m(u)\leqslant 2||u||_m^2+C(1+||u||_0^{4m}).\]
   Denote $C_m=C(1+||u(0)||_0^{4m})+C||f||_m^2$, then from (\ref{(A.1)}) and above, we deduce
   \[\frac{d}{dt}(\mathcal{J}_m(u)-C_m)\leqslant \frac{1}{2}\epsilon (\mathcal{J}_m(u)-C_m),\]
   thus \[\mathcal{J}_m(u)-C_m\leqslant e^{\frac{1}{2}\epsilon t}[\mathcal{J}_m(u(0))-C_m],\]
    so\[||u(t)||_m^2\leqslant 4||u(0)||_m^2 e^{\frac{1}{2}\epsilon t}+C_m.\]
    This prove Proposition  0.3.\quad $ \square$
    
    \section*{Appendix }
    Consider the following system of ordinary differential equations:
    \[\dot{x}=Y(x),\quad x(0)=x_0\in \mathbb{R}^n,\]
    where   $Y(x)=(Y_1(x),\cdots,Y_n(x)): \mathbb{R}^n\to\mathbb{R}^n$ is a continuously differentiable map. Let $F(t,x)$ be a (local) flow  determined by this equation.
    \smallskip
    
    \noindent {\bf Theorem A} (Liouville).  Let $B(x_1,\cdots,x_n)$ be a continuous differentiable function on  $\mathbb{R}^n$. For the Borel measure $d\mu=B(x)dx$ in $\mathbb{R}^n$
   and any bounded open set $A\subset \mathbb{R}^n$, we have 
    \[\frac{d}{dt}\mu(F(t,A))=\int_{F(t,A)}\Big[\sum_{i=1}^n\frac{\partial (B(x)Y_i(x))}{\partial x_i}\Big]dx, \quad t\in (-T,T),\]
    where $T>0$ is such that $F(t,x)$ is well defined and bounded for any $t\in~(-T,T)$ and $x\in A$.
    \smallskip
    
    For $B=const$ this result is well known.     For its proof for  a non-constant density $B$ see e.g. \cite{15, 16}.
    
    \section*{Acknowlegement}
        First of all, I want to   thank  my PhD supervisor S. Kuksin for formulation of the problem and  guidance. I am grateful  to  A. Boritchev for useful suggestions.  Finally, I would like  to thank  all of the staff and faculty at  CMLS of Ecole Polytechnique   for help.
        
 \section*{Reference}
\bibliography{pKdV_HG.bib}
\end{document}